\documentclass[nocompress]{spie}

\usepackage{graphicx}
\usepackage{subcaption}
\usepackage{multirow}
\usepackage{caption}
\usepackage{subcaption}
\usepackage{hyperref}
\usepackage{amsmath,amsfonts,amssymb}
\usepackage[table]{xcolor}
\usepackage{import}
\usepackage{float}
\usepackage{wrapfig}
\usepackage[export]{adjustbox}
\usepackage{caption}
\usepackage{verbatim}
\usepackage{tabularx}
\usepackage{makecell}

\begin{document}
\title{Reinforcement learning for graph theory \\ Parallelizing Wagner’s approach}
\author[a]{Alix Bouffard}
\author[a]{Jane Breen}
\affil[a]{Ontario Tech University}

\date{}
\maketitle

\section*{Abstract}
     Our work applies reinforcement learning to construct counterexamples concerning conjectured bounds on the spectral radius of the Laplacian matrix of a graph. We expand upon the re-implementation of Wagnar’s approach by Stevanovic et al. with the ability to train numerous unique models simultaneously and a novel redefining of the action space to adjust the influence of the current local optimum on the learning process. \cite{Wagner} \cite{Stev1}
    

\section{Introduction}
    \subsection{Contributions}
    A shortcoming of the cross-entropy method is its tendency to become stuck in local minima. Our implementation attempts to address this issue by executing models in parallel, effectively performing the function of multiple local searches. To this end, we have re-implemented and modularized Stevonovic's code-base from scratch. Another avenue of possible improvement mentioned by the previous authors pertains to the rate of convergence during the training process. We hope to develop a more exploratory approach by starting a portion of each batch at the nearest local maxima for each individual model. To enable this, the action space of the model has been changed.
    
    \subsection{Table of Contents}
    This thesis is segmented into chapters, hyper-linked here for the readers convenience. The introduction is comprised of two parts, \hyperref[sec:BKG]{Background Knowledge} and \hyperref[sec:PW]{Previous Works}. Background knowledge provides a brief overview of the basics of graph theory, the Laplacian matrix, and spectral graph theory. Previous works reviews the history of computer-generated conjectures and counterexamples, along with past implementations of the code-base this thesis used as a reference. The \hyperref[sec:method]{Methodology} section begins with a review of reinforcement learning and the cross-entropy optimization method as it pertains to graph theory. It then addresses the redefined action space and the implementation of threaded models. \hyperref[sec:Results]{Results} details the findings of this thesis and provides direction for future research. \hyperref[Appendix]{Appendix A} lists the counter-examples found by this thesis as well as the corresponding graphs and adjacency matrices.

\newpage
    \subsection{Background Knowledge}
    \label{sec:BKG}

    \begin{wrapfigure}{r}{0.55\textwidth}
        \includegraphics[width=\linewidth]{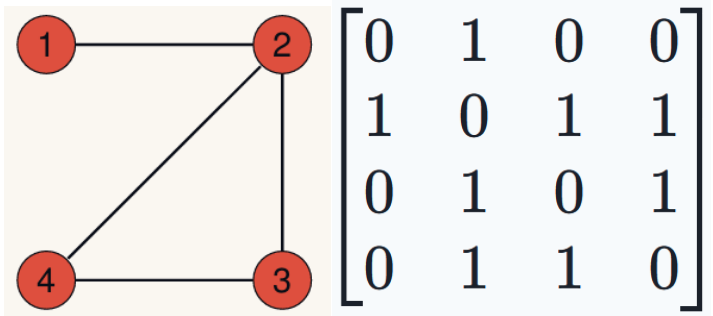}
        \caption{A four vertex graph and its adjacency matrix}
        \label{fig:graph_adj}
    \end{wrapfigure}
    
    A graph G is a mathematical object defined in terms of a vertex set V(G), $\{1, 2, 3\}$, and an edge set E(G) consisting of pairs of vertices, $\{\{1, 2\}, \{1, 3\}, \{2, 3\}\}$. Vertices u and v are said to be adjacent to one another if $\{u, v\} \in E(G)$. The degree of a vertex is the number of vertices adjacent to it. The standard way of representing a graph is using a matrix called the adjacency matrix. Let G be a graph with vertices labelled 1, 2, … $n$. The adjacency matrix of G is defined to be the $n$×$n$ matrix A such that the ($i$, $j$) entry of A is 1 if there is an edge between $i$ and $j$, and 0 otherwise. Note that the adjacency matrix is symmetric.
    \newline
    \newline
    \noindent
    Graphs can be used to represent many real-world systems such as the internet, transportation networks, food webs, power grids, and neural networks. These networks can have thousands of vertices and tens of thousands of edges, which makes it difficult to extract useful information about the graph through either real-world observation or computer simulation. Spectral graph theory is a field of study which focuses on the relationship between the structural properties of a graph and the set of eigenvalues of a matrix associated with that graph. From the eigenvalues of the adjacency matrix one can deduce the lower bound on the maximum degree of a node in the graph, the maximum distance between any two nodes, and the number of cycles of length three, among many other characteristics.
    \newline
    \newline
    Using the adjacency matrix, the Laplacian matrix of a graph can be formed. The eigenvalues of the Laplacian matrix can be used to extract further information about a graph. A few uses of the eigenvalues of the Laplacian include determining the number of connected components of a graph and partitioning communities of nodes. \cite{Autograffix}
    \vspace{1em}
    \begin{wrapfigure}{l}{0.7\textwidth}
        \centering
        \includegraphics[width=\linewidth]{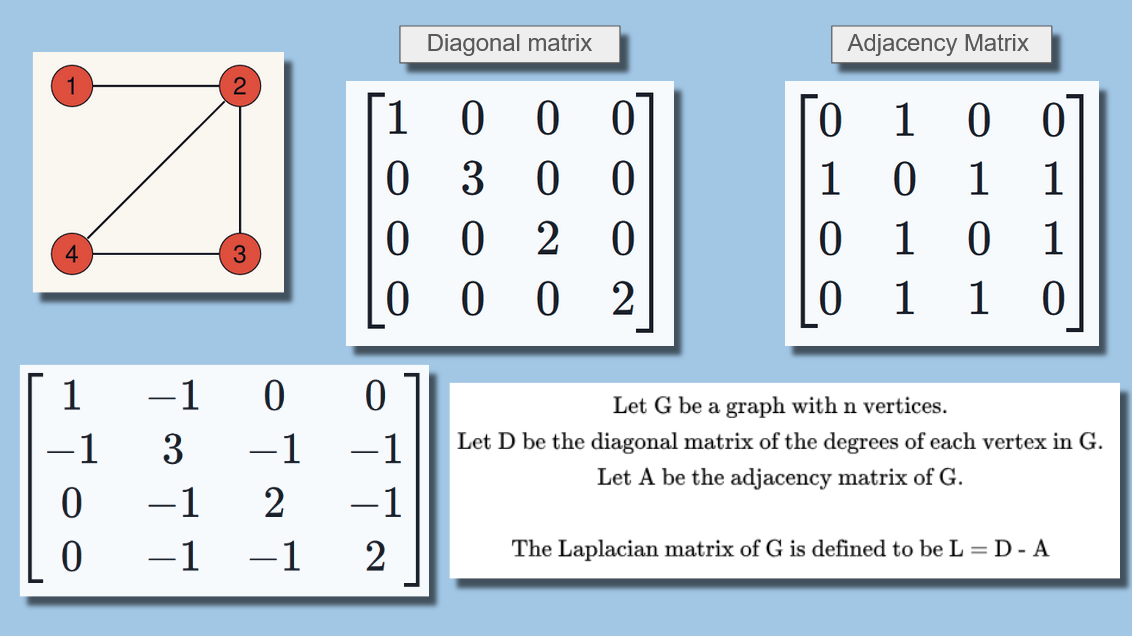}
        \caption{The Laplacian matrix is displayed in the lower left}
        \label{fig:Laplacian}
    \end{wrapfigure}
    
    The Laplacian matrix is positive semi-definite, meaning all of its eigenvalues are real and nonnegative. In particular, the largest eigenvalue µ of the Laplacian matrix, also called the spectral radius of L, has found numerous applications so far. It has been used to determine the first ionization potential of alkanes, in combinatorial optimization to provide an upper bound on the size of the maximum cut in a graph, and in communication networks to provide a lower bound on the edge-forwarding index.\cite{brankov2006automated} \cite{faught2024nordhaus}
    \newline
    \newpage
    \subsection{Previous Works}
    \label{sec:PW}
    Wagner developed an initial approach for applying a reinforcement learning cross-entropy method to generate counterexamples for several conjectures in extremal combinatorics and graph theory.\cite{Wagner} This method was notable as it drastically reduced the computational requirement of evaluating larger graphs. As there are over 11 million connected graphs on 10 vertices, it is infeasible to exhaustively search every graph as some conjectures have counterexamples in the hundreds of vertices. In the context of graph theory, this method successfully addressed multiple conjectures, including disproving the preservation of transmission regularity under co-spectrality of the distance Laplacian, clarifying the relationship between proximity and distance eigenvalues, and examining conjectures involving the sum of the largest eigenvalue and the graph's matching number.  These conjectures are particularly well-suited to a reinforcement learning approach as they naturally yield, or can be manipulated to yield, negative rewards until a counterexample is discovered. 
\begin{figure}[H]
        \centering
        \includegraphics[width=0.4\linewidth]{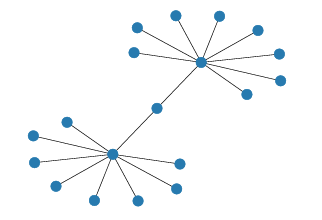}
        \caption{A graph on 19 vertices satisfying $\lambda_1 + \mu < \sqrt{n - 1} + 1$}
        \label{fig:Wagner example}
    \end{figure}

    For the conjectures on transmission regularity and the sum of eigenvalue and matching number, the algorithm directly produced valid counterexamples. Although the model did not yield a direct counterexample to the proximity conjecture, the graphs that achieved the maximum reward gave a clear indication on the probable structure that a counterexample might take. Wagner was then able to explore similarly structured graphs until he found a counter example. Thus, even when a direct solution is not obtained, this method can provide significant value in guiding theoretical exploration.
    \begin{figure}[H]
        \begin{subfigure}{0.4\textwidth}
            \includegraphics[width=\textwidth]{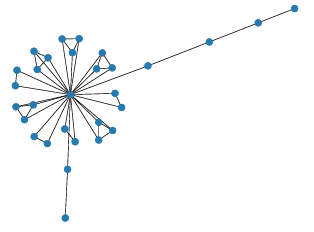} 
            \caption{The top performing graph by Wagner's implementation}
            \label{fig:wag2-1}
        \end{subfigure}
        \hfill
        \begin{subfigure}{0.6\textwidth}
            \includegraphics[width=\textwidth]{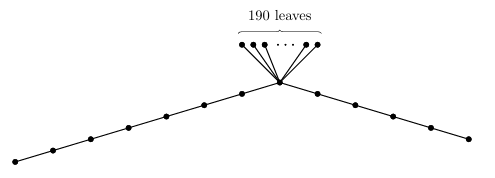}
            \caption{Corresponding counter-example}
            \label{fig:wag2-2}
        \end{subfigure}
        
        \caption{Counter-example for the conjecture: $\pi + \delta_{\lfloor\frac{2D}{3}\rfloor} > 0$}
        \label{fig:Wagner example 2}
    \end{figure}
\newpage
    The conjecture concerning the sum of the largest eigenvalue and the matching number was originally generated by AutoGraphiX, a software designed to automatically generate conjectures regarding relationships between graph invariants.\cite{Autograffix}\cite{Autocon} However, there has been much work on automated and computer-aided generation of conjectures and as such there are many other programs.\cite{TxGraffiti} Brankov et al observed a recurring structure in upper bounds for the Laplacian spectral radius: these bounds are often expressed as maxima of functions involving a vertex's degree ($d_i$) and the average degree of its neighbours  ($m_i$). They developed a program that generated conjectures structured as functions relating these parameters. For example: $\rho(L) - \text{max}(d_i + \sqrt{d_im_i} )\leq 0$
    \newline
    \newline
    From the thousands of conjectures generated, 190 vertex-maximum and 297 edge-maximum conjectures were assessed on all connected graphs with up to 9 vertices, and no counterexample were found. In the appendix, they listed 68 open conjectured bounds.\cite{brankov2006automated} It is those bounds for which Stevanovic et al. were driven to refine Wagner’s approach.
    \newline
    \newline
    The focus of the work of Stevanovic et al. was to enhance the readability, stability, and computational efficiency of Wagner’s codebase.\cite{Stev1}\cite{Stev2} A drawback of Wagner’s original implementation was that it required specific combinations of Python and TensorFlow versions. To avoid this, Stevanovic re-implemented Wagner’s approach from scratch. Their revised design separates the reinforcement learning agent from the reward computation, enabling future researchers to easily adapt the code to their use case by specifying their reward logic in a standalone Python file. Although graph-specific reward calculations can be performed using conventional Python libraries such as NetworkX and NumPy, they used graph6java to increase performance, a dedicated Java source file imported into the Python environment.\cite{graph6java} This hybrid approach provided significant runtime advantages, accelerating performance by approximately 3 to 5 times over Wagner.

\newpage

\section{Methodology}
\vspace{1em}
\label{sec:method}
    \noindent
    \begin{minipage}[t]{0.4\textwidth}
        \vspace*{-1\baselineskip}
        \centering
        \vspace{3.5em}
        \includegraphics[width=\linewidth]{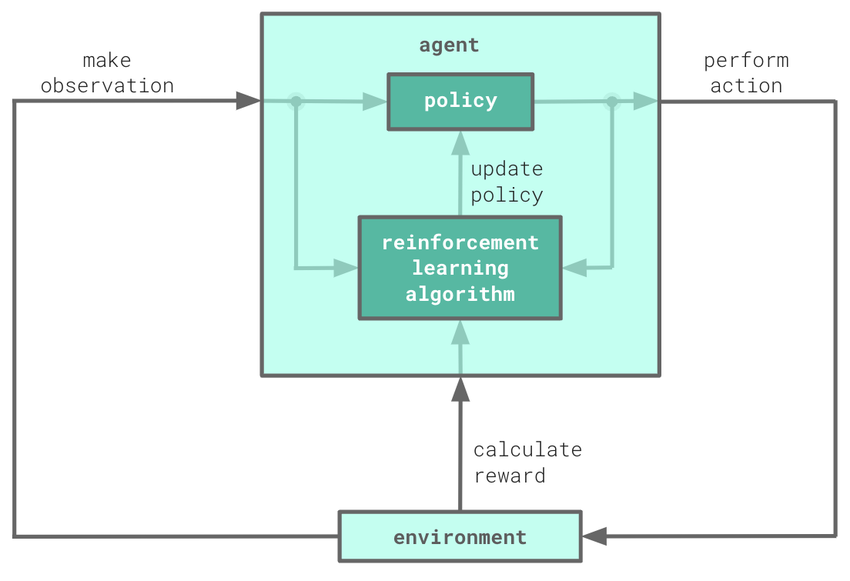}
        \captionof{figure}{A reinforcement learning algorithm enacts a policy on an environment, makes observations, and optimizes itself according to a reward.\cite{Fowler}}
    \end{minipage}
    \hfill
    \begin{minipage}[t]{0.55\textwidth}
        \raggedright
        Cross-entropy is a basic evolutionary algorithm, where populations are generated, evaluated, and iterated upon.\cite{cetutorial}  Populations are generated from  the policy network, which is updated based on the characteristics of the top performing individuals in the previous generation. Optionally, the top performing individuals can be carried forward to the next generation without modification to maintain the quality of the generation. It is a simple reinforcement learning method, used for non-complex environments where a single final reward given at the terminal state is thought to be sufficient to evaluate the performance of the model.
        \newline
        \newline
        The model attempts to maximize the reward received at the end of each generation. As the conjectures are given as an
        \begin{wrapfigure}{r}{0.5\textwidth}
            \vspace{-1em}
            \centering
            \includegraphics[width=\linewidth]{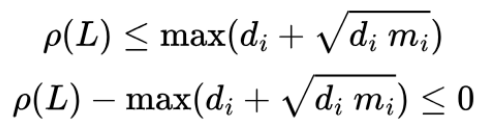}
            \caption{Reward function}
        \end{wrapfigure}
        upper bound on the spectral radius, the reward is chosen to be the spectral radius minus the conjecture.
        \vspace{-2em}
        \newline
        \newline
        \newline
        This has the benefit of being a negative reward until a counterexample is found.
    \end{minipage}
    \newline
    \newline
    \noindent
    In our implementation, the top 10\% performing graphs are passed to the model for learning and the top 5\% are duplicated and appended to the next generation. It is important to note that the highest rewarded graph is carried through generations in perpetuity until displaced and that the surviving individuals are also included in the learning segment.

    \subsection{Model environment}

    \noindent
    \begin{minipage}[t]{0.5\textwidth}
        \vspace*{-1\baselineskip}
        \includegraphics[width=0.45\textwidth]{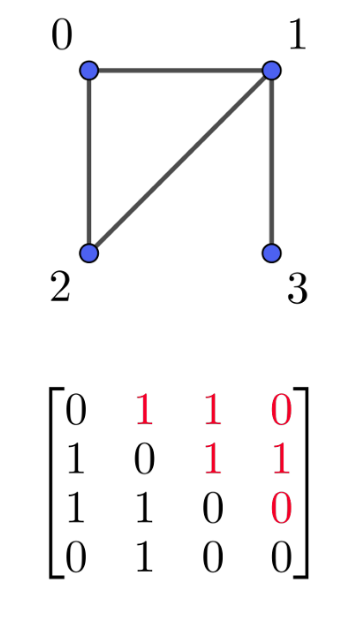}
        \hfill
        \includegraphics[width=0.45\textwidth]{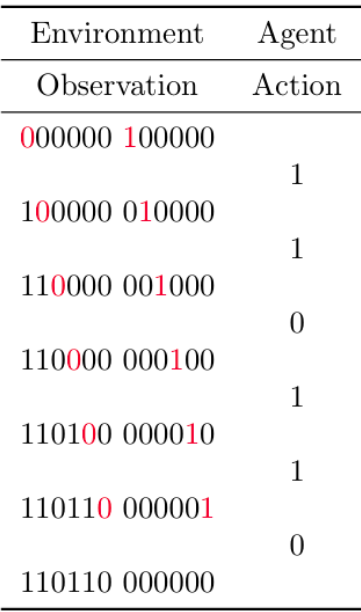}
    
        \captionof{figure}{Left: A graph and its adjacency matrix. The section of the adjacency matrix used as the model environment is highlighted in red. Right: The adjacency matrix being constructed by the model.}
        \label{fig:Wagner example 2}
    \end{minipage}
    \hfill
    \begin{minipage}[t]{0.47\textwidth}
        The model constructs a graph by returning a boolean action based on the observation matrix it receives at each step. The observation matrix consists of the current elements of the adjacency matrix above the diagonal, listed in row-wise order, followed by a one-hot encoding of the current vertex the model is making a decision on. Previous work used an intuitive action space of inserting the chosen action directly to the adjacency matrix. This was enabled by having the model build each graph from the zero matrix at each generation.
    \end{minipage}

    \enlargethispage{3\baselineskip}
    \begin{wrapfigure}{r}{0.25\textwidth}
            \vspace*{-10\baselineskip}
            \centering
            \includegraphics[width=\linewidth]{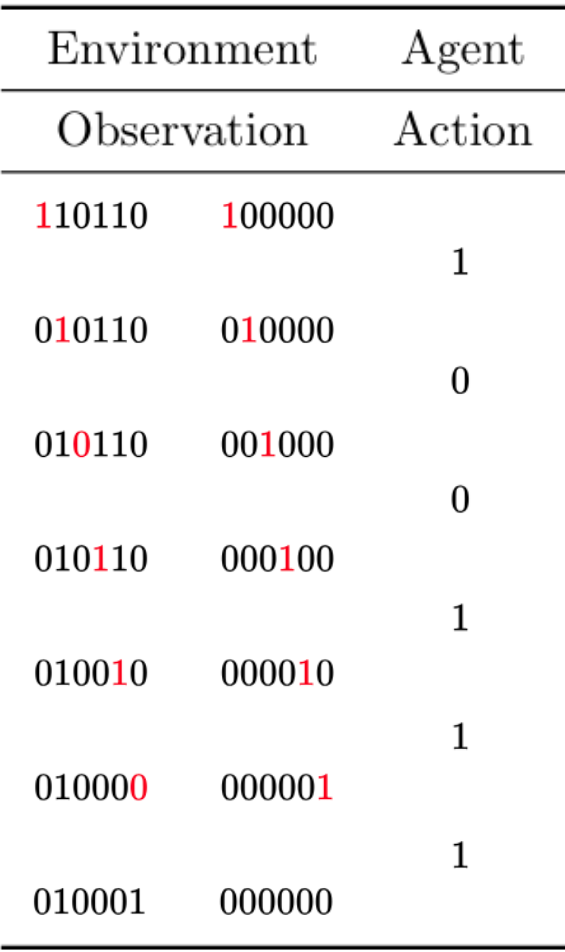}
    \end{wrapfigure}
    
    In an effort to encourage the model to produce graphs of similar structure to the previous top performing graphs, we initialize a portion of each generation with the previous generation's top observation matrix. This required the action space to be redefined as the choice to alter the current observation. This was implemented as a bit-wise exclusive or between the model's action and the observation matrix. In the case of the zero-initialized observation matrix, the redefined action space is congruent to previous works.
    \clearpage

    \newpage
    \subsection{Parallel processing}
    The cross-entropy method displays common drawbacks found in the field of reinforcement learning: lower efficiency as the dimensionality of solution space increases, and the penchant to get stuck in local optima. \cite{DecentCEM} The work of Zhang et al. attempts to remedy these problems by decentralizing the training process by running multiple models in parallel. Each model runs independently from one another, dispersed along the solution space through the initial randomization of each network.

    \noindent
    \begin{minipage}[t]{0.45\textwidth}
    As illustrated in the inlaid figure, the centralized sampling distribution exhibits a bias toward the sub-optimal solutions near top right, due to the global top-k ranking. \cite{DecentCEM} In comparison, a decentralized approach could maintain diversity due to its local sampling instances. Implementing this technique required a complete restructuring of Stevanovic's code base. This restructuring enabled the modularization of executing model actions, evaluating rewards, and updating model policies. 
    \end{minipage}
    \hfill
    \begin{minipage}[t]{0.5\textwidth}
        \vspace*{-1\baselineskip}
        \includegraphics[width=0.45\textwidth]{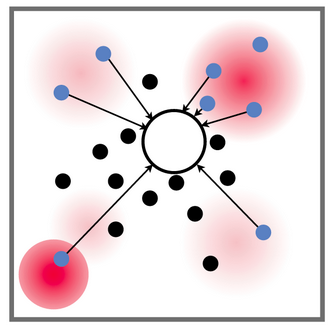}
        \hfill
        \includegraphics[width=0.44\textwidth]{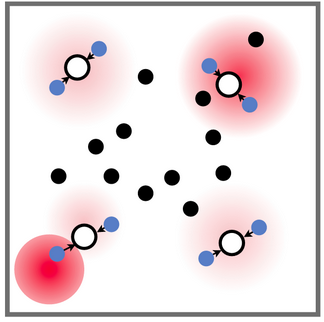}
        \captionof{figure}{Left: Centralized and De-centralized models.}
    \end{minipage}
    \noindent
    \newline
    \newline
    These changes reduced the time to process a generation by up to a third in comparison to Stevanovic's implementation. To facilitate a comparison between a single model implementation, the batch size of the single model is divided among the instances.

    \subsection{Network architecture}
    Each model is a sequential, fully connected network from the pytorch library. The size of the input layer is dictated by the length of the observation matrix. The default architecture for the hidden network are two layers of 72 and 12 nodes. The output layer is a pair of nodes whose normalized outputs are used directly as probabilities for the model actions. The default learning rate is 0.002 used in an ADAM optimizer with a cross-entropy loss objective.\cite{ADAM} The activation function of each node is GELU, an alteration from previous implementations.\cite{GELU} Each threaded model is able to have its own unique hyper-parameters, though this function was not utilized.

      \begin{figure}[H]
            \centering
            \includegraphics[width=0.9\textwidth]{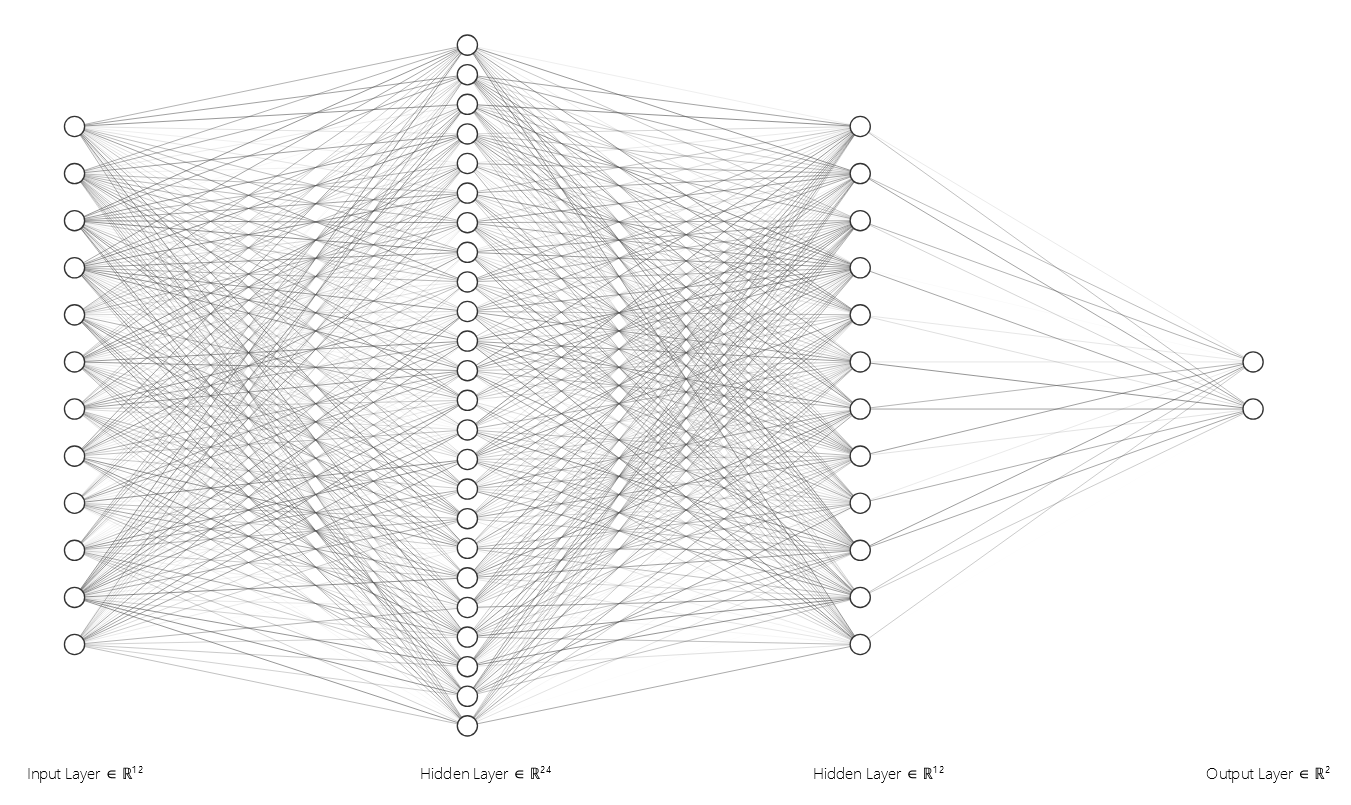}
            \caption{Example network for a 4-node graph. 12 node observation/input layer. 24 and 12 node hidden layers. 2 node output layer.}
    \end{figure}

    \newpage{}
    \section{Results}
    \label{sec:Results}
    
    \begin{figure}[H]
        \centering
        \includegraphics[width=\linewidth]{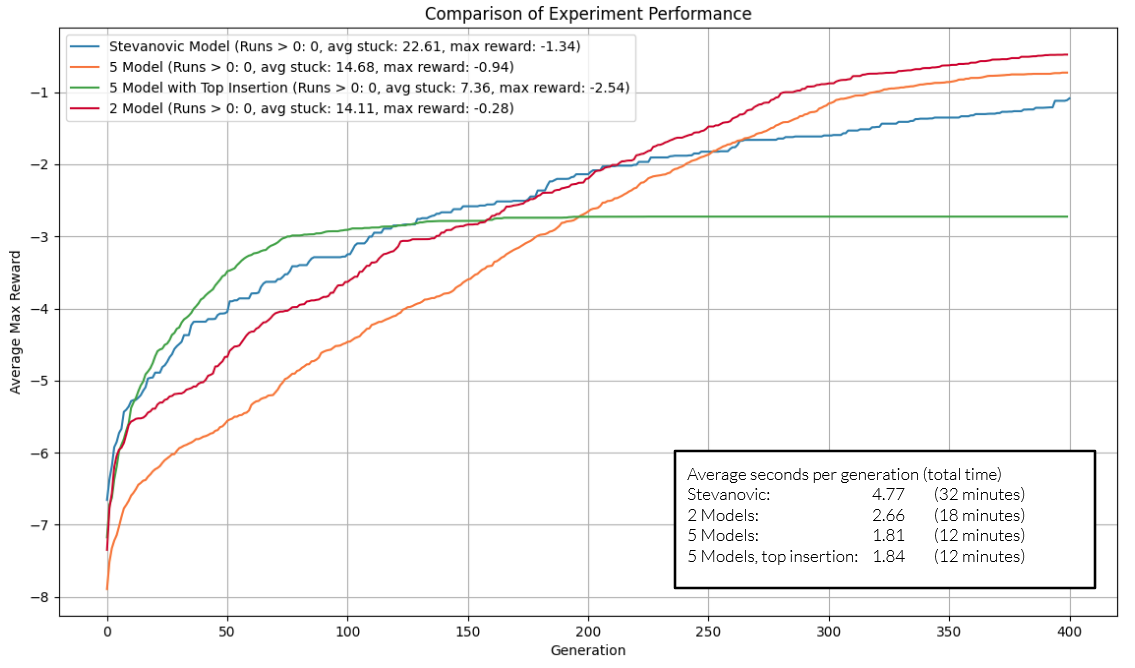}
        \caption{Comparison of various implementation performances on conjecture 3}
    \end{figure}
    
    These experiments are the average of 15 runs of 400 generations each with a total batch size of 1000. In addition to finding new counterexamples, our model shows promise in scaling. The tests were performed for conjectured bound 3 on 20 nodes. The bound has a \hyperref[sec:graph 2]{counterexample} which was found during a prior test using a batch size of 200 on a single model. Though the counterexample was not found during these experimental tests, it is clear that the proposed method has a higher average reward and executes faster. Both the 5 model and 2 model implementations outperformed the original implementation. A possible reason for the initial under-performance of the 5 model is the batch size. With 5 models, the batch size is 200 compared to 500, this limits the breadth of graphs able to be generated in a generation. Initializing a quarter of the batch with the current top iteration provides an initial boost but quickly leads to getting stuck in a local optima. Future research might explore starting with a high insertion percent and tapering off as the generations progress.

    \begin{figure}[H]
        \centering
        \includegraphics[width=\linewidth]{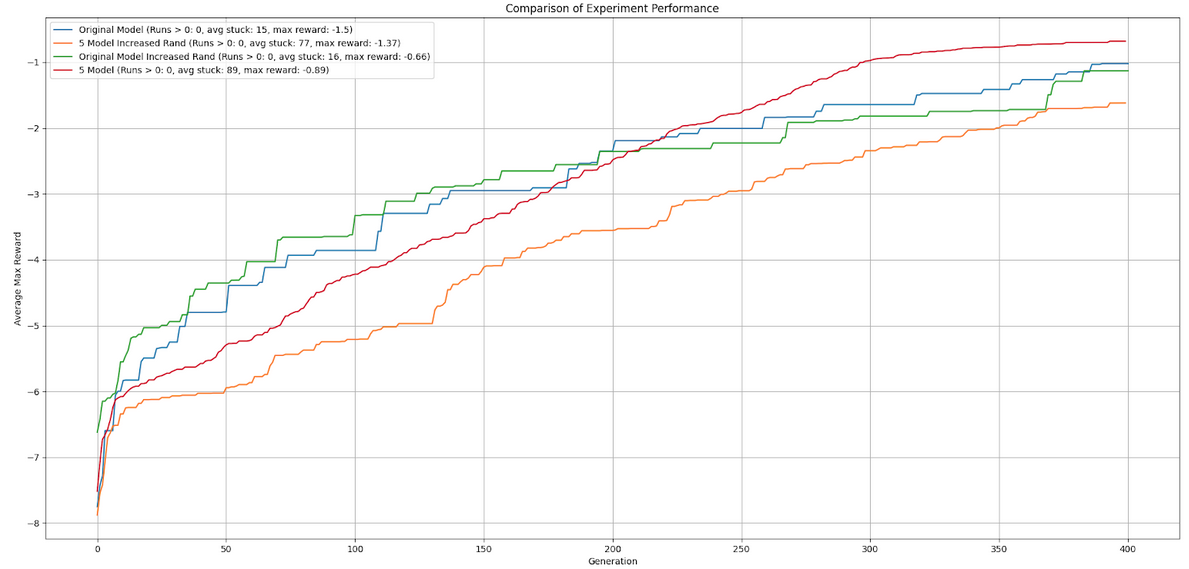}
        \caption{Comparison of the impact of randomness on performances on conjecture 3}
    \end{figure}
       
    In Stevanovic's implementation, the number of graphs randomized in any given step starts at a base level.\cite{Stev1} If the model is stuck at a maxima for a set number of generations, the amount of randomness is increased. In our implementation, the number of random actions taken are constant throughout the generations, regardless of the time spent at the current maxima. One of the variables we adjusted was increasing the number of random graphs per batch from 0.05\% to 1\%. The goal of this change was to increase the rate of improvement near the beginning of the run as well as preventing any prolonged stretches at a local maxima by encouraging exploration. This slight increase in randomization led to a marked decrease in performance in our implementation, and a minor decrease in performance in Stevanovic's. It is possible that during the modularization of the randomization function the effect of threaded models on the RNG seed was overlooked, leading to unexpected results.
       
    \newpage

       {\rowcolors{3}{gray!25}{white!}
        \begin{tabular}{ |p{2.5cm}|p{2.5cm}|p{2.5cm}|p{2.5cm}|p{2.5cm}|  }
        \hline
        \multicolumn{5}{|c|}{Counter examples} \\
        \hline
        Conjectured Bound & Graph 2 & Graph 66 & Graph 41 & Graph 65 \\
        \hline
        2 & $\circ$     & \_         & \_         & \_         \\
        3 & \checkmark  & \checkmark & \_         & \_         \\
        15 & \checkmark & \checkmark & \_         & \_         \\
        28 & \checkmark & \checkmark & \_         & \_         \\
        29 & \checkmark & \checkmark & \_         & \_         \\
        31 & \checkmark & \checkmark & \_         & \_         \\
        32 & $\circ$    & $\circ$    & \_         & \_         \\
        36 & \checkmark & \checkmark & \_         & \_         \\
        41 & \_         & \_         & \checkmark & \_         \\
        43 & \checkmark & \checkmark & \checkmark & \_         \\
        49 & \checkmark & \checkmark & \checkmark & \_         \\
        51 & \_         & \_         & \checkmark & \_         \\
        52 & \checkmark & \checkmark & \checkmark & \_         \\
        53 & \checkmark & \checkmark & \checkmark & \_         \\
        54 & \checkmark & \checkmark & \checkmark & \_         \\
        55 & \checkmark & \checkmark & \checkmark & \_         \\
        57 & \checkmark & \checkmark & \checkmark & \_         \\
        58 & \checkmark & \checkmark & \checkmark & \_         \\
        59 & \checkmark & \checkmark & \_         & \_         \\
        60 & \checkmark & \checkmark & \_         & \_         \\
        61 & $\circ$    & $\circ$    & \_         & \_         \\
        62 & \checkmark & \checkmark & \_         & \_         \\
        63 & \checkmark & \checkmark & \_         & \_         \\
        64 & \checkmark & \checkmark & \_         & \_         \\
        65 & \_         & \_         & \_         & \checkmark \\
        66 & \_         & \checkmark & \_         & \_         \\
        67 & \checkmark & \_         & \_         & \_         \\
        68 & \_         & \_         & \_         & \checkmark  \\
        \hline
        \end{tabular}}

    \vspace{1em}
    The work of Stevanovic provided a unique graph for each counter-example, this thesis condenses these results to four graphs. Each graph is named for the first unique counter-example found compared to the previous graphs. The $\checkmark$ symbol indicates a counter-example to conjectures that were previously closed, while $\circ$ denotes new counter-examples. The graphs and corresponding conjectures can be found in \hyperref[Appendix]{Appendix A}.
       
\section{Concluding remarks}
    In conclusion, our model outperforms the original implementation in both time-complexity and the number of counterexamples found. Potential future research can focus on restricting the model to certain classes of graphs and an exploration of the number of isomorphic graphs generated in each generation. The code-base is publicly available at \url{https://github.com/Alix-B/Parallelizing-Wagners-Approach}.

\section*{Acknowledgements}
    \centering
    This thesis is submitted in partial fulfillment of the requirements for the degree of \\ Bachelor of Science in Applied and Industrial Mathematics, Faculty of Science.

\newpage
\nocite{Annealing}
\nocite{charton2024patternboost}
\nocite{ramsey}
\nocite{kanervisto2020action}
\nocite{jaderberg2017population}
\bibliography{bibliography}
\bibliographystyle{spiebib}

\newpage
\appendix
\section{Counter-examples}
\label{Appendix}
    Following the diagram and adjacency matrix of each graph, the conjectures to which they are a counter-example to are provided. That is, the spectral radius of the Laplacian matrix of the graph is larger than the stated conjecture. The variables of the conjectures are as follows. Given a graph $G$, $V$ is its vertex set. $d_v$ is the degree of a vertex $v$, $m_v$ is the average degrees of the neighbours of $v$. For conjectures taken over the maximum of all pairs of adjacent vertices $v$ and $j$, $d_j$ and $m_j$ are used.

    \newpage
    \subsection{Graph 2}
    \label{sec:graph 2}
    \begin{figure}[H]
    \centering
    \begin{minipage}[t]{0.3\linewidth}
        \centering
        \includegraphics[width=\linewidth]{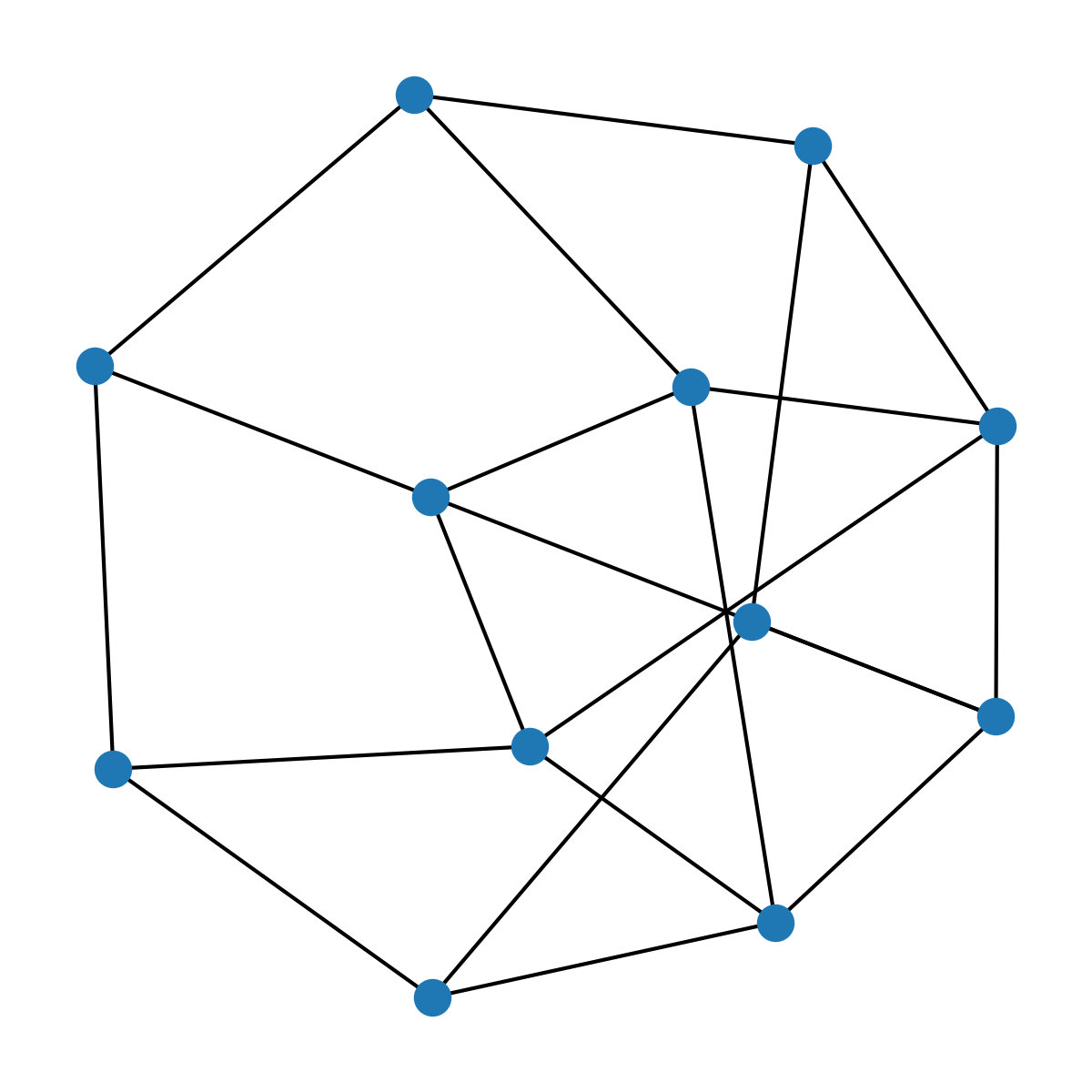}
    \end{minipage}
    \hfill
    \begin{minipage}[t]{0.6\linewidth}
        \centering
        \vspace{-15em}
        \begin{verbatim}
        [[0 0 0 1 1 0 0 0 1 0 0 0]
         [0 0 1 1 0 0 0 0 1 0 1 0]
         [0 1 0 0 0 1 0 1 0 0 0 1]
         [1 1 0 0 0 0 0 1 0 0 0 1]
         [1 0 0 0 0 0 1 1 0 0 0 0]
         [0 0 1 0 0 0 0 0 1 1 0 0]
         [0 0 0 0 1 0 0 0 0 1 1 0]
         [0 0 1 1 1 0 0 0 0 0 1 0]
         [1 1 0 0 0 1 0 0 0 0 0 0]
         [0 0 0 0 0 1 1 0 0 0 0 1]
         [0 1 0 0 0 0 1 1 0 0 0 1]
         [0 0 1 1 0 0 0 0 0 1 1 0]]
        \end{verbatim}
    \end{minipage}
    \end{figure}

    \begin{tabularx}{\textwidth}{@{}X X@{}}
        \makecell[l]{\textbf{Conjecture 2:}\\[0.3ex] $\displaystyle \max_{v \in V} \frac{2m_v^2}{d_v}$} &
        \makecell[l]{\textbf{Conjecture 3:}\\[0.3ex] $\displaystyle \max_{v \in V} \frac{m_v^2}{d_v} + m_v$} \\
        \makecell[l]{\textbf{Conjecture 15:}\\[0.3ex] $\displaystyle \max_{v \in V} \sqrt{\frac{4m_v^3}{d_v}}$} \\
        \makecell[l]{\textbf{Conjecture 28:}\\[0.3ex] $\displaystyle \max_{v \in V} \sqrt{\frac{4m_v^4}{d_v^2} + 2d_vm_v}$} &
        \makecell[l]{\textbf{Conjecture 29:}\\[0.3ex] $\displaystyle \max_{v \in V} \sqrt{m_v^2 + \frac{3m_v^3}{d_v}}$} \\
        \makecell[l]{\textbf{Conjecture 31:}\\[0.3ex] $\displaystyle \max_{v \in V} \frac{4m_v^2}{m_v+d_v}$} &
        \makecell[l]{\textbf{Conjecture 32:}\\[0.3ex] $\displaystyle \max_{v \in V} \frac{\sqrt{m_v^3(m_v+3d_v)}}{d_v}$} \\
        \makecell[l]{\textbf{Conjecture 36:}\\[0.3ex] $\displaystyle \max_{v \sim j} \frac{2(m_v^2+m_j^2)}{d_v+d_j}$} &
        \makecell[l]{\textbf{Conjecture 43:}\\[0.3ex] $\displaystyle \max_{v \sim j} 2 + \sqrt{3(m_v^2+m_j^2)-2m_vm_j-4(d_v+d_j)+4}$} \\
        \makecell[l]{\textbf{Conjecture 49:}\\[0.3ex] $\displaystyle \max_{v \sim j} 2 + \sqrt{2(m_v^2+m_j^2)+(d_v-d_j)^2-4(d_v+d_j)+4}$} &
        \makecell[l]{\textbf{Conjecture 52:}\\[0.3ex] $\displaystyle \max_{v \sim j} 2 + \sqrt{\sqrt{8(m_v^4+m_j^4)-8(d_v^2+d_j^2)+4}-4(d_v+d_j)+6}$} \\
        \makecell[l]{\textbf{Conjecture 53:}\\[0.3ex] $\displaystyle \max_{v \sim j} 2 + \sqrt{\sqrt{8(m_v^4+m_j^4)-8(d_vm_v+d_jm_j)+4}-4(d_v+d_j)+6}$} \\
        \makecell[l]{\textbf{Conjecture 54:}\\[0.3ex] $\displaystyle\max_{v \sim j} 2 +\sqrt{2(m_v^2+m_j^2)+(d_vm_v+d_jm_j)-(d_v^2+d_j^2)-4(d_v+d_j)+4}$} \\
        \makecell[l]{\textbf{Conjecture 55:}\\[0.3ex] $\displaystyle\max_{v \sim j} 2 +\sqrt{3(m_v^2+m_j^2)-(d_v^2+d_j^2)-4(m_v+m_j)+4}$} &
        \makecell[l]{\textbf{Conjecture 57}\\[0.3ex] $\displaystyle \max_{v \sim j} 2 +\sqrt{2(m_v^2+m_j^2)-8\frac{d_v^2+d_j^2}{m_v+m_j}+4}$} \\
        \makecell[l]{\textbf{Conjecture 58:}\\[0.3ex] $\displaystyle \max_{v \sim j} 2 +\sqrt{2(m_v^2+m_vm_j+m_j^2)-(d_vm_v+d_jm_j) -4(d_v+d_j)+4}$} \\
        \makecell[l]{\textbf{Conjecture 59:}\\[0.3ex] $\displaystyle \max_{v \sim j} \frac{2(m_v^2+m_vm_j+m_j^2)-(d_v^2+d_j^2)}{m_v+m_j}$} &
        \makecell[l]{\textbf{Conjecture 60:}\\[0.3ex] $\displaystyle \max_{v \sim j} 2 +\sqrt{2(m_v^2+m_vm_j+m_j^2)-(d_v^2+d_j^2) -4(d_v+d_j)+4}$} \\
        \makecell[l]{\textbf{Conjecture 61}\\[0.3ex] $\displaystyle \max_{v \sim j} \frac{2(m_v^2+m_j^2)}{2+\sqrt{2((d_v-1)^2+(d_j-1)^2)}}$} &
        \makecell[l]{\textbf{Conjecture 62:}\\[0.3ex] $\displaystyle \max_{v \sim j} 2 +\sqrt{m_v^2+4m_vm_j+m_j^2-2d_vd_j -4(d_v+d_j)+4}$} \\
        \makecell[l]{\textbf{Conjecture 63}\\[0.3ex] $\displaystyle \max_{v \sim j} d_v + d_j + m_v + m_j - 4\frac{d_vd_j}{m_v+m_j}$} &
        \makecell[l]{\textbf{Conjecture 64:}\\[0.3ex] $\displaystyle \max_{v \sim j} \frac{m_vm_j(d_v + d_j)}{d_vd_j}$} \\
        \makecell[l]{\textbf{Conjecture 67:}\\[0.3ex] $\displaystyle \max_{v \sim j} \frac{(m_v+m_j)(d_vm_v + d_jm_j)}{2d_vd_j}$}
    \end{tabularx}

    \newpage
    \enlargethispage{3\baselineskip}
    \subsection{Graph 66}
    \begin{minipage}[t]{0.3\linewidth}
        \centering
        \includegraphics[width=\linewidth]{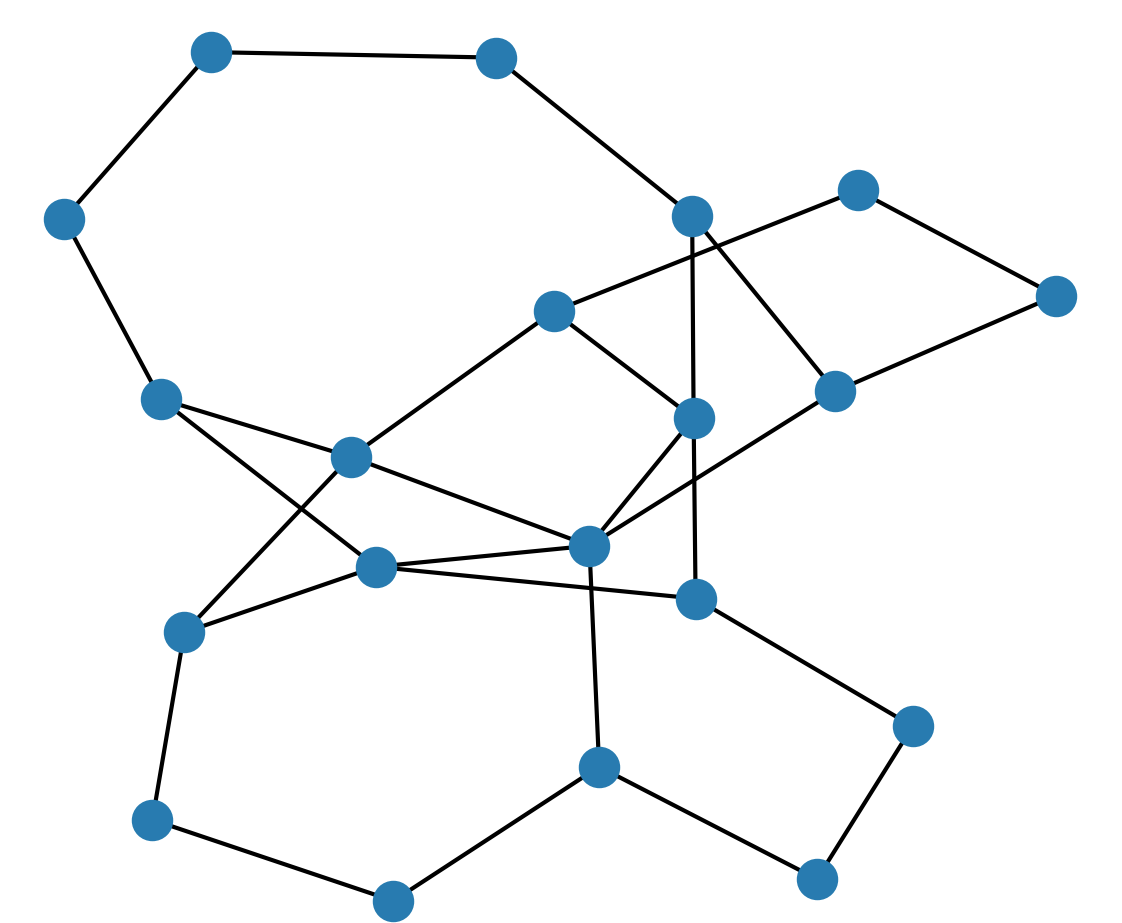}
    \end{minipage}
    \hfill
    \begin{minipage}[t]{0.6\linewidth}
        \centering
        \vspace{-12em}
        \begin{verbatim}
        [[0 0 1 0 1 0 0 1 0 0 0 1 0 0 0 0 0 0 0 0]
         [0 0 1 0 0 0 0 1 1 0 0 0 0 0 0 0 1 0 0 0]
         [1 1 0 0 0 0 0 0 0 0 1 0 0 0 0 0 0 0 0 0]
         [0 0 0 0 0 0 0 0 0 0 0 1 0 0 0 0 0 0 0 1]
         [1 0 0 0 0 1 0 0 0 0 0 0 1 0 0 0 0 0 0 0]
         [0 0 0 0 1 0 0 0 0 1 0 0 0 0 0 0 0 0 0 0]
         [0 0 0 0 0 0 0 0 0 0 0 0 0 1 0 0 0 0 1 0]
         [1 1 0 0 0 0 0 0 0 0 0 0 1 0 0 1 0 0 1 0]
         [0 1 0 0 0 0 0 0 0 0 0 0 0 1 0 1 0 0 0 0]
         [0 0 0 0 0 1 0 0 0 0 0 0 0 0 1 0 0 0 0 0]
         [0 0 1 0 0 0 0 0 0 0 0 0 0 0 0 0 0 1 0 0]
         [1 0 0 1 0 0 0 0 0 0 0 0 0 0 0 1 0 0 0 0]
         [0 0 0 0 1 0 0 1 0 0 0 0 0 0 0 0 0 0 0 1]
         [0 0 0 0 0 0 1 0 1 0 0 0 0 0 0 0 0 0 0 0]
         [0 0 0 0 0 0 0 0 0 1 0 0 0 0 0 0 1 0 0 0]
         [0 0 0 0 0 0 0 1 1 0 0 1 0 0 0 0 1 0 0 0]
         [0 1 0 0 0 0 0 0 0 0 0 0 0 0 1 1 0 0 0 0]
         [0 0 0 0 0 0 0 0 0 0 1 0 0 0 0 0 0 0 1 0]
         [0 0 0 0 0 0 1 1 0 0 0 0 0 0 0 0 0 1 0 0]
         [0 0 0 1 0 0 0 0 0 0 0 0 1 0 0 0 0 0 0 0]]
        \end{verbatim}
    \end{minipage}

    \begin{tabularx}{\textwidth}{@{}X X@{}}
        \makecell[l]{\textbf{Conjecture 3:}\\[0.3ex] $\displaystyle \max_{v \in V} \frac{m_v^2}{d_v} + m_v$} &
        \makecell[l]{\textbf{Conjecture 15:}\\[0.3ex] $\displaystyle \max_{v \in V} \sqrt{\frac{4m_v^3}{d_v}}$} \\
        \makecell[l]{\textbf{Conjecture 28:}\\[0.3ex] $\displaystyle \max_{v \in V} \sqrt{\frac{4m_v^4}{d_v^2} + 2d_vm_v}$} &
        \makecell[l]{\textbf{Conjecture 29:}\\[0.3ex] $\displaystyle \max_{v \in V} \sqrt{m_v^2 + \frac{3m_v^3}{d_v}}$} \\
        \makecell[l]{\textbf{Conjecture 31:}\\[0.3ex] $\displaystyle \max_{v \in V} \frac{4m_v^2}{m_v+d_v}$} &
        \makecell[l]{\textbf{Conjecture 32:}\\[0.3ex] $\displaystyle \max_{v \in V} \frac{\sqrt{m_v^3(m_v+3d_v)}}{d_v}$} \\
        \makecell[l]{\textbf{Conjecture 36:}\\[0.3ex] $\displaystyle \max_{v \sim j} \frac{2(m_v^2+m_j^2)}{d_v+d_j}$} &
        \makecell[l]{\textbf{Conjecture 43:}\\[0.3ex] $\displaystyle \max_{v \sim j} 2 + \sqrt{3(m_v^2+m_j^2)-2m_vm_j-4(d_v+d_j)+4}$} \\
        \makecell[l]{\textbf{Conjecture 49:}\\[0.3ex] $\displaystyle \max_{v \sim j} 2 + \sqrt{2(m_v^2+m_j^2)+(d_v-d_j)^2-4(d_v+d_j)+4}$} &
        \makecell[l]{\textbf{Conjecture 52:}\\[0.3ex] $\displaystyle \max_{v \sim j} 2 + \sqrt{\sqrt{8(m_v^4+m_j^4)-8(d_v^2+d_j^2)+4}-4(d_v+d_j)+6}$} \\
        \makecell[l]{\textbf{Conjecture 53:}\\[0.3ex] $\displaystyle \max_{v \sim j} 2 + \sqrt{\sqrt{8(m_v^4+m_j^4)-8(d_vm_v+d_jm_j)+4}-4(d_v+d_j)+6}$} \\
        \makecell[l]{\textbf{Conjecture 54:}\\[0.3ex] $\displaystyle\max_{v \sim j} 2 +\sqrt{2(m_v^2+m_j^2)+(d_vm_v+d_jm_j)-(d_v^2+d_j^2)-4(d_v+d_j)+4}$} \\
        \makecell[l]{\textbf{Conjecture 55:}\\[0.3ex] $\displaystyle\max_{v \sim j} 2 +\sqrt{3(m_v^2+m_j^2)-(d_v^2+d_j^2)-4(m_v+m_j)+4}$} &
        \makecell[l]{\textbf{Conjecture 57}\\[0.3ex] $\displaystyle \max_{v \sim j} 2 +\sqrt{2(m_v^2+m_j^2)-8\frac{d_v^2+d_j^2}{m_v+m_j}+4}$} \\
        \makecell[l]{\textbf{Conjecture 58:}\\[0.3ex] $\displaystyle \max_{v \sim j} 2 +\sqrt{2(m_v^2+m_vm_j+m_j^2)-(d_vm_v+d_jm_j) -4(d_v+d_j)+4}$} \\
        \makecell[l]{\textbf{Conjecture 59:}\\[0.3ex] $\displaystyle \max_{v \sim j} \frac{2(m_v^2+m_vm_j+m_j^2)-(d_v^2+d_j^2)}{m_v+m_j}$} &
        \makecell[l]{\textbf{Conjecture 60:}\\[0.3ex] $\displaystyle \max_{v \sim j} 2 +\sqrt{2(m_v^2+m_vm_j+m_j^2)-(d_v^2+d_j^2) -4(d_v+d_j)+4}$} \\
        \makecell[l]{\textbf{Conjecture 61}\\[0.3ex] $\displaystyle \max_{v \sim j} \frac{2(m_v^2+m_j^2)}{2+\sqrt{2((d_v-1)^2+(d_j-1)^2)}}$} &
        \makecell[l]{\textbf{Conjecture 62:}\\[0.3ex] $\displaystyle \max_{v \sim j} 2 +\sqrt{m_v^2+4m_vm_j+m_j^2-2d_vd_j -4(d_v+d_j)+4}$} \\
        \makecell[l]{\textbf{Conjecture 63}\\[0.3ex] $\displaystyle \max_{v \sim j} d_v + d_j + m_v + m_j - 4\frac{d_vd_j}{m_v+m_j}$} &
        \makecell[l]{\textbf{Conjecture 64:}\\[0.3ex] $\displaystyle \max_{v \sim j} \frac{m_vm_j(d_v + d_j)}{d_vd_j}$} \\
        \makecell[l]{\textbf{Conjecture 66:}\\[0.3ex] $\displaystyle \max_{v \sim j} \frac{m_v^2 +4m_vm_j+m_j^2-(d_vm_v+d_jm_j)}{d_v+d_j}$}
    \end{tabularx}
    \clearpage

    \newpage
    \enlargethispage{3\baselineskip}
    \subsection{Graph 41}
    \begin{figure}[H]
    \centering
    \begin{minipage}[t]{0.3\linewidth}
        \centering
        \includegraphics[width=\linewidth]{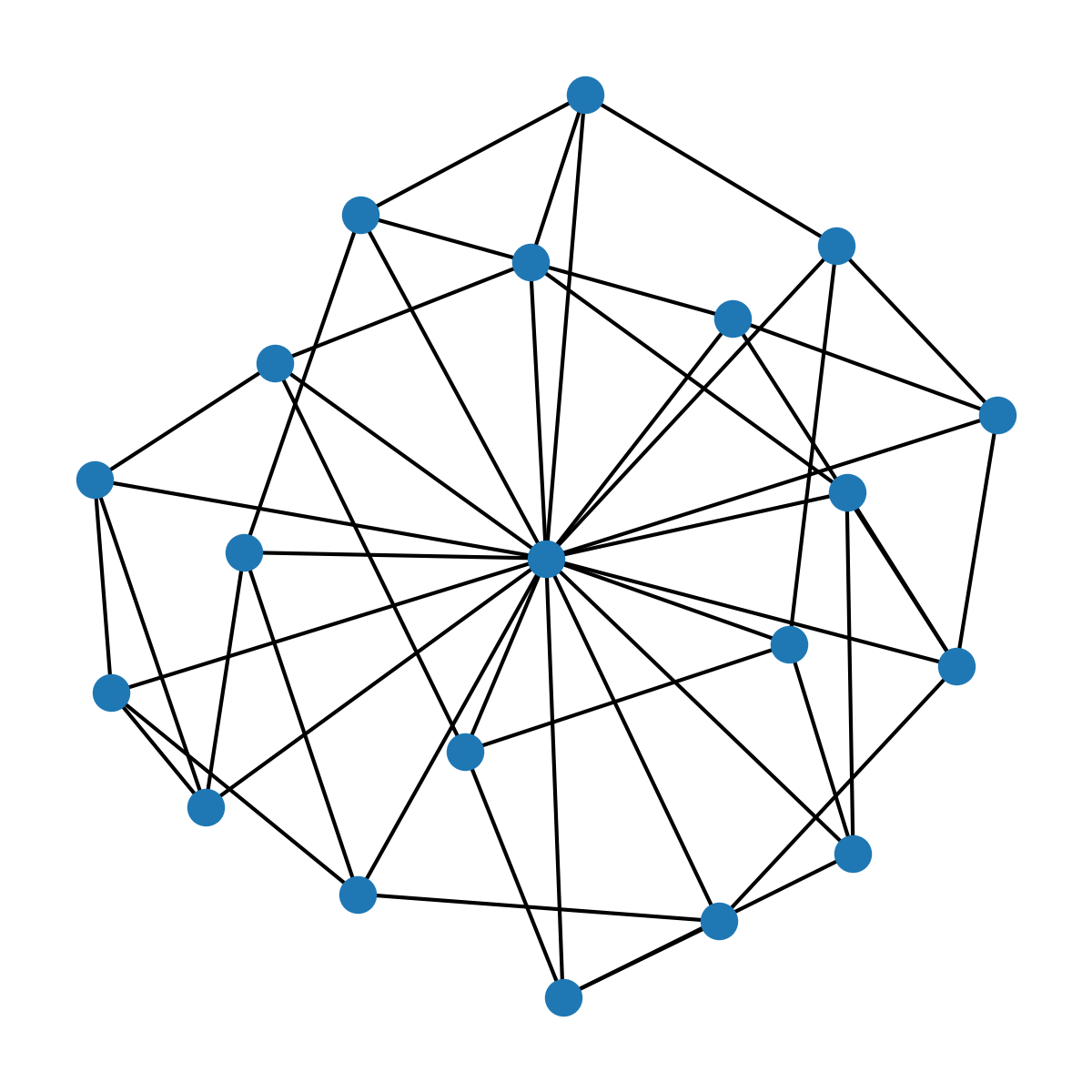}
    \end{minipage}
    \hfill
    \begin{minipage}[t]{0.6\linewidth}
        \centering
        \vspace{-15em}
        \begin{verbatim}
        [[0 1 1 1 1 1 1 1 1 1 1 1 1 1 1 1 1 1 1 1]
         [1 0 0 1 0 0 0 1 0 0 0 0 0 1 0 0 0 0 0 0]
         [1 0 0 0 0 0 1 0 0 0 0 0 1 0 1 0 0 0 0 0]
         [1 1 0 0 0 1 0 0 0 0 0 0 0 0 1 0 0 0 0 0]
         [1 0 0 0 0 1 0 0 0 0 0 0 0 0 0 0 1 0 1 0]
         [1 0 0 1 1 0 0 0 0 0 0 0 1 0 0 0 0 0 0 0]
         [1 0 1 0 0 0 0 0 0 0 0 0 0 0 1 0 1 1 0 0]
         [1 1 0 0 0 0 0 0 1 0 0 0 0 0 0 0 0 0 1 0]
         [1 0 0 0 0 0 0 1 0 0 0 0 0 1 0 0 0 1 0 0]
         [1 0 0 0 0 0 0 0 0 0 0 0 1 0 0 1 0 0 0 1]
         [1 0 0 0 0 0 0 0 0 0 0 1 0 0 0 1 0 0 0 1]
         [1 0 0 0 0 0 0 0 0 0 1 0 0 0 0 1 0 0 1 0]
         [1 0 1 0 0 1 0 0 0 1 0 0 0 0 0 0 0 0 0 0]
         [1 1 0 0 0 0 0 0 1 0 0 0 0 0 0 0 1 0 0 0]
         [1 0 1 1 0 0 1 0 0 0 0 0 0 0 0 0 0 0 0 0]
         [1 0 0 0 0 0 0 0 0 1 1 1 0 0 0 0 0 0 0 0]
         [1 0 0 0 1 0 1 0 0 0 0 0 0 1 0 0 0 0 0 0]
         [1 0 0 0 0 0 1 0 1 0 0 0 0 0 0 0 0 0 0 1]
         [1 0 0 0 1 0 0 1 0 0 0 1 0 0 0 0 0 0 0 0]
         [1 0 0 0 0 0 0 0 0 1 1 0 0 0 0 0 0 1 0 0]]
        \end{verbatim}
    \end{minipage}
    \end{figure}

    \begin{tabularx}{\textwidth}{@{}X X@{}}
        \makecell[l]{\textbf{Conjecture 41:}\\[0.3ex] $\displaystyle \max_{v \sim j} 2 + (m_v+m_j)-(d_v+d_j)+\sqrt{2(d_v^2+d_j^2)-4(m_v+m_j)+4}$} \\
        \makecell[l]{\textbf{Conjecture 43:}\\[0.3ex] $\displaystyle \max_{v \sim j} 2 + \sqrt{3(m_v^2+m_j^2)-2m_vm_j-4(d_v+d_j)+4}$} \\
        \makecell[l]{\textbf{Conjecture 49:}\\[0.3ex] $\displaystyle \max_{v \sim j} 2 + \sqrt{2(m_v^2+m_j^2)+(d_v-d_j)^2-4(d_v+d_j)+4}$} \\
        \makecell[l]{\textbf{Conjecture 51:}\\[0.3ex] $\displaystyle \max_{v \sim j} 2(m_v+m_j)-4\frac{m_vm_j}{d_v+d_j}$} \\
        \makecell[l]{\textbf{Conjecture 52:}\\[0.3ex] $\displaystyle \max_{v \sim j} 2 + \sqrt{\sqrt{8(m_v^4+m_j^4)-8(d_v^2+d_j^2)+4}-4(d_v+d_j)+6}$} \\
        \makecell[l]{\textbf{Conjecture 53:}\\[0.3ex] $\displaystyle \max_{v \sim j} 2 + \sqrt{\sqrt{8(m_v^4+m_j^4)-8(d_vm_v+d_jm_j)+4}-4(d_v+d_j)+6}$} \\
        \makecell[l]{\textbf{Conjecture 54:}\\[0.3ex] $\displaystyle\max_{v \sim j} 2 +\sqrt{2(m_v^2+m_j^2)+(d_vm_v+d_jm_j)-(d_v^2+d_j^2)-4(d_v+d_j)+4}$} \\
        \makecell[l]{\textbf{Conjecture 55:}\\[0.3ex] $\displaystyle\max_{v \sim j} 2 +\sqrt{3(m_v^2+m_j^2)-(d_v^2+d_j^2)-4(m_v+m_j)+4}$} \\
        \makecell[l]{\textbf{Conjecture 57}\\[0.3ex] $\displaystyle \max_{v \sim j} 2 +\sqrt{2(m_v^2+m_j^2)-8\frac{d_v^2+d_j^2}{m_v+m_j}+4}$} \\
        \makecell[l]{\textbf{Conjecture 58:}\\[0.3ex] $\displaystyle \max_{v \sim j} 2 +\sqrt{2(m_v^2+m_vm_j+m_j^2)-(d_vm_v+d_jm_j) -4(d_v+d_j)+4}$}
    \end{tabularx}
    \clearpage

    \newpage
    \subsection{Graph 65}
    \begin{figure}[H]
    \centering
    \begin{minipage}[t]{0.3\linewidth}
        \centering
        \includegraphics[width=\linewidth]{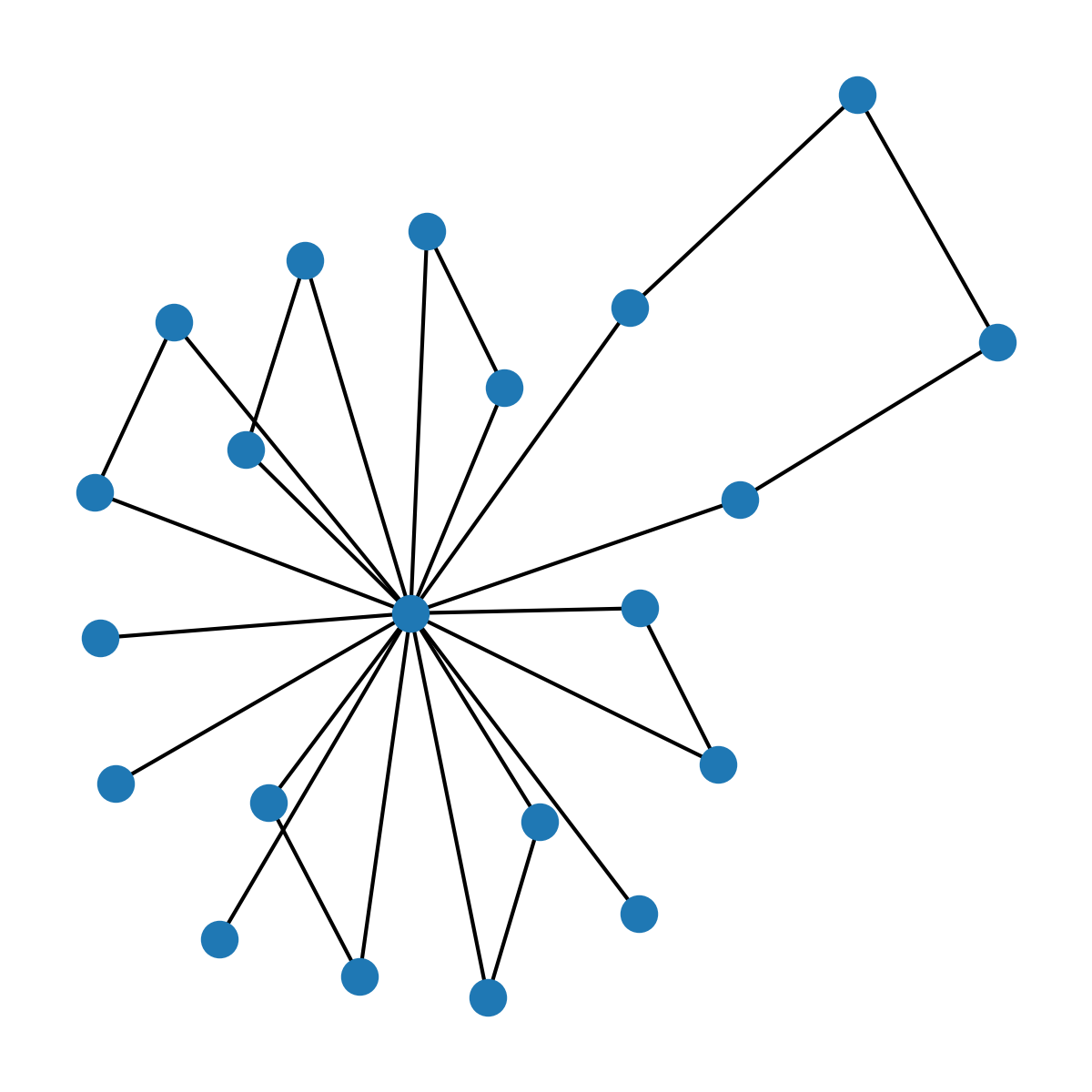}
    \end{minipage}
    \hfill
    \begin{minipage}[t]{0.6\linewidth}
        \centering
        \vspace{-15em}
        \begin{verbatim}
        [[0 0 0 1 1 1 1 1 1 1 1 1 1 1 1 1 1 1 1 1 1]
         [0 0 1 1 0 0 0 0 0 0 0 0 0 0 0 0 0 0 0 0 0]
         [0 1 0 0 0 0 0 1 0 0 0 0 0 0 0 0 0 0 0 0 0]
         [1 1 0 0 0 0 0 0 0 0 0 0 0 0 0 0 0 0 0 0 0]
         [1 0 0 0 0 0 0 0 0 0 0 0 0 1 0 0 0 0 0 0 0]
         [1 0 0 0 0 0 1 0 0 0 0 0 0 0 0 0 0 0 0 0 0]
         [1 0 0 0 0 1 0 0 0 0 0 0 0 0 0 0 0 0 0 0 0]
         [1 0 1 0 0 0 0 0 0 0 0 0 0 0 0 0 0 0 0 0 0]
         [1 0 0 0 0 0 0 0 0 0 0 0 0 0 0 0 0 1 0 0 0]
         [1 0 0 0 0 0 0 0 0 0 0 0 0 0 0 0 0 0 1 0 0]
         [1 0 0 0 0 0 0 0 0 0 0 0 0 0 0 0 0 0 0 0 0]
         [1 0 0 0 0 0 0 0 0 0 0 0 0 0 0 0 1 0 0 0 0]
         [1 0 0 0 0 0 0 0 0 0 0 0 0 0 0 0 0 0 0 0 0]
         [1 0 0 0 1 0 0 0 0 0 0 0 0 0 0 0 0 0 0 0 0]
         [1 0 0 0 0 0 0 0 0 0 0 0 0 0 0 0 0 0 0 1 0]
         [1 0 0 0 0 0 0 0 0 0 0 0 0 0 0 0 0 0 0 0 0]
         [1 0 0 0 0 0 0 0 0 0 0 1 0 0 0 0 0 0 0 0 0]
         [1 0 0 0 0 0 0 0 1 0 0 0 0 0 0 0 0 0 0 0 0]
         [1 0 0 0 0 0 0 0 0 1 0 0 0 0 0 0 0 0 0 0 0]
         [1 0 0 0 0 0 0 0 0 0 0 0 0 0 1 0 0 0 0 0 0]
         [1 0 0 0 0 0 0 0 0 0 0 0 0 0 0 0 0 0 0 0 0]]
        \end{verbatim}
    \end{minipage}
    \end{figure}

    \noindent
    \begin{tabularx}{\textwidth}{@{}X X@{}}
        \makecell[l]{\textbf{Conjecture 65:}\\[0.3ex] $\displaystyle \max_{v \sim j} \frac{(m_v+m_j)(d_vm_v+d_jm_j)}{2m_vm_j}$} &
        \makecell[l]{\textbf{Conjecture 68:}\\[0.3ex] $\displaystyle \max_{v \sim j} 2 + \sqrt{(m_v-m_j)^2 +4d_vd_j-4(m_v+m_j)+4}$}
    \end{tabularx}

\end{document}